\documentclass[12pt,leqno]{article}
\usepackage{amsmath,amsthm,amsfonts,amssymb,amscd}
\usepackage[latin1]{inputenc}
\usepackage{graphics}

\begin {document}

\title{\bf Example of a foliation for Ueda type one}
\date{04/08/2023}
\author{Paulo Sad}

\maketitle

\begin {abstract} 
\noindent We give an example of a projective 2-dimensional foliation which is regular along a curve of genus 3 and has  Ueda type  one with normal bundle   of order two. 

\end {abstract}

\noindent {\bf Introduction}

\vspace {0.1in}

In this Introduction we present a short list of results of the so called Ueda Theory, developed in \cite {U} with focus in local properties. Later on this matter was studied in \cite {N} from a global point of view, and recently appered again in the work of \cite{P}.

Let $C$ be a smooth, compact, connected holomorphic curve contained in some holomorphic surface $M$ satisfying  $C.C=0$ ($C.C$ stands for the self intersection number of $C$ inside $M$). We consider a type $\nu$ coordinate system $\{\psi_i: \mathbb D \times \mathbb D \rightarrow M\}$ for a neighborhood of $C$:  
 
\begin {description} 

\item $z_k= \psi_k^{-1}\circ \psi_i (z_i,y_i)$ \\
                                                        
\item $y_k={\lambda}_{ik}(z_i)y_i + f_{ik}(z_i)y_i^{\nu+1}+\dots $

\end {description}

\noindent for some $\nu \ge 1$.

\vspace {0.1in}
The normal bundle $N_C$ to $C$ is given by

\begin {description}
 
\item                                                  $z_k= \psi_k^{-1}\circ \psi_i (z_i,0)$ \\
\item                                                  $y_k={\lambda}_{ik}(z_i)y_i$  

\end {description}

\vspace {0.1in}
\noindent Since $C.C=0$,  we may assume that $\lambda_{ik}(z_i)=\lambda_{ik} \in \mathbb C^*$.

\vspace {0.1in}
The collection $\{f_{ik}\}$ is a cocycle of the sheaf of holomorphic functions of $N_C^{-1}$; when it is a coboundary, we may replace the type $\nu$ coordinate system by a $\nu+1$ coordinate system. The curve $C$ has finite type $n\in \mathbb N$ when there exists a system of coordinates of type $n$ whose associated cocycle is not a coboundary; we say that $n$ is the {\it Ueda type} of $C$. The Ueda type of $C$ is {\it infinite} when the cocycle associated to any coordinate system is a coboundary.

\vspace {0.1in}
The finite type is characterized by the following result:

\vspace {0.1in}
\noindent {\bf Theorem 1}\,\, If $C$ has finite type then it has a fundamental system of pseudo-concave neighborhoods.

\vspace{0.1in}
\noindent In particular, there are no compact curves in a neighborhood of $C$ disjoint of $C$.

When the Ueda type is infinite, the extra condition of $N_C$ being a torsion bundle leads to the following

\vspace{0.1in}
\noindent {\bf Theorem 2}\,\, If $C$ has infinite type and its normal bundle $N_C$ is a torsion bundle of order $k \in \mathbb N$ then there exists a holomorphic foliation by compact
leaves in a neighborhood of $C$ such that $C$ is a leaf (in fact the foiation can be seen as a fibration which has $C$ as a multiple fiber of multiplicity $k$).

\vspace {0.1in}
Both  results appeared in \cite {U}. We state now the results of \cite {N} and \cite {P} respectively in the case $M$ is a projective surface.

\vspace {0.1in}
\noindent {\bf Theorem 3}\, \, Assume that $N_C$ is a torsion bundle of order $k\in \mathbb N$ and let $n\in \mathbb N \cup \{\infty\}$ be the Ueda type of $C$. If $n > k$ then in fact $n=\infty$ (and Theorem 2 applies).

\vspace{0.1in}
\noindent {\bf Theorem 4}\,\, Assume $N_C$ is a torsion bundle of order $k \in \mathbb N$ and let $n \in \mathbb N$ be the Ueda type of $C$. If $n=k$ then there exists a foliation in $M$ which has $C$ as a leaf.

\vspace {0.1in}
We remind that $C$ being a leaf means in particular that it contains no singularity of the foliation.

\vspace {0.1in} 
Theorem 4 leads to the following question: suppose $N_C$ is a torsion bundle of order $k\in \mathbb N$ and let $n\in \mathbb N$ be the Ueda type of $C$;  what can be said about the existence of a foliation having $C$ as a leaf if $n<k$? We present here an example in the positive direction when $n=1$ and $k=2$ (correcting a statement of \cite{Sa}). 
 
\vspace{0.1in}
We have in mind the following picture: we start with the projective plane $\mathbb P^2$ with some foliation which has a smooth algebraic invaariant curve $S$ of degree $d$; this curve $S$ contains singularities $p_1\dots,p_l$ with tangent multiplicities (of the foliation along $S$) $m_1,\dots,m_l$. According to \cite{S}, the degree $0$ divisor $\sum_{j=1}^lm_jp_j-LD_{\infty}$ is a principal divisor ($d.L=\sum_{j=1}^lm_j$ and $D_{\infty}$ is the divisor obtained intersecting $S$ with some line, which we fix at infinity). Since a generic choice of points in $S$ produces no "resonances" there is no foliation having such points as singularities in the invariant curve $S$. In particular, when we blow-up $\mathbb P^2$ at $d^2$ points in $S$ chosen in a generic way the strict transform $\widehat S$ of $S$ can not be a leaf of a foliation in the transformed surface $\widehat {\mathbb P^2}$.

\vspace {0.1in}
In this paper we will work with a particular Neeman's example. This family of examples (see \cite {N}) is constructed in the following way: we start with a fixed quartic $S$ and blow it up at 16 points which lie at two different bitangent lines (4 points $P_1,\dots,P_4$) and some cubic (12 points $P_5,\dots,P_{16}$). The Ueda type of $\widehat S$ is 1, and the normal bundle $N_{\widehat S}$ is associated to the divisor $\sum_{j=1}^{16}P_j-4D_{\infty}$. We see that $2(\sum_{j=1}^{16}P_j-4D_{\infty}$) is a principal divisor, so that $N_{\widehat S}$ is a torsion bundle of order 2. Such a resonance is realized by a foliation in $\mathbb P^2$ of degree 10; it has 16 singularities of tangent multiplicities equal to 2 along $S$. The group of holonomy of $S\setminus \{P_1,\dots,P_{16}\}$ (which is the same of $\widehat S$) is infinite, although the local holonomy maps at the singularities are all igual to the identity map.

\vspace {0.2in}
\noindent {\bf Example}

\vspace {0.2in}

Let $F(X,Y)$ be a degree four polynomial which defines a smooth projective quartic $S$ transverse to the line at infinity; we assume also that $\{F=0\}\cap\{F_X=0\}$ has 12 different points and $l_1(X,Y)=Y+X=0$ and $l_2(X,Y)=Y+X+b=0$ are bitangent lines. We put

\vspace{0.1in}
\begin{description}
\item $G=l_1l_2F_X^2$
\item $h=\dfrac{(l_1l_2)_XF_Y}{2}+\dfrac{(l_1l_2)_YF_X}{2}+l_1l_2F_{XY}=\dfrac{(l_1l_2)_X(F_X+F_Y)}{2}+l_1l_2F_{XY}$
\item $k=-F_X-\dfrac{(l_1l_2)_XF_{XY}}{2}$
\item $h^{\prime}=-(l_1l_2)_XF_X-l_1l_2F_{XX}$
\item $k^{\prime}= F_X+\dfrac{(l_1l_2)_XF_{XX}}{2}$ 
\item $A=hF_X+kF$
\item $B=h^{\prime}F_X+k^{\prime}F$
\end{description}

\vspace{0.1in}
The foliation $\cal F$ we will study is given by

$$
GdF+F(BdX-AdY)=0
$$

\noindent Let $\widehat {\cal F}$   be the transformed foliation from $\cal F$  and $\widehat S$ be the transformed curve from $S$ after blowing-up at the points of $\{G=0\}\cap S$.

\vspace {0.2in}
\noindent {\bf Theorem.} $\widehat S$ is a leaf of $\widehat {\cal F}$.

\vspace{0.1in} First of all we need to compute the degree of the foliation $\cal F$.

\vspace {0.1in}
\noindent {\bf Lemma.}  $  \cal F$  has degree 10.

\vspace {0.1in}
\noindent Proof: it ammounts to showing that the homogeneous polynomial of degree 12 in the development of

\begin{description}
\item $X(GF_X+FB)+Y(GF_Y-FA)=$
\item $XGF_X+YGF_Y+F[X(h^{\prime}F_X+k^{\prime}F)-Y(hF_X+kF)]$
\end {description}

\noindent is in fact equal to 0. In what follows, all the equalities refer to homogeneous polynomials of top degree. For example, $XGF_X+YGF_Y=4GF$.

\vspace {0.1in}

Let us compute $X(h^{\prime}F_X+k^{\prime }F)-Y(hF_X+kF)$. We have

\begin{description}

\item $X(h^{\prime}F_X+k^{\prime}F)-Y(hF_X+kF)=$

\item $=-X(l_1l_2F_X)_XF_X-Y\left[\dfrac{(l_1l_2)_XF_Y}{2}+\dfrac{(l_1l_2)_YF_X}{2}+l_1l_2F_{XY}\right]F_X$

\item $+X\left[\dfrac{(l_1l_2)_{XX} F_X}{2}+ \dfrac{(l_1l_2)_X F_{XX}}{2}\right]F
-Y\left[-F_X-\dfrac{(l_1l_2)_X F_{XY}}{2}\right]F$

\item $=-X(l_1l_2F_X)_XF_X-Y\left[\dfrac{(l_1l_2)_XF_Y}{2}+\dfrac{(l_1l_2)_YF_X}{2}+l_1l_2F_{XY}\right]F_X$

\item $+(X+Y)FF_X+ \left[\dfrac{X(l_1l_2)_X F_{XX}}{2}+\dfrac{Y(l_1l_2)_XF_{XY}}{2}\right]F$

\item $=-X(l_1l_2F_X)_XF_X-Y\left[\dfrac{(l_1l_2)_XF_Y}{2}+\dfrac{(l_1l_2)_YF_X}{2}+l_1l_2F_{XY}\right]F_X$

\item $+(X+Y)FF_X +\dfrac{3(l_1l_2)_XF_XF}{2}$

\item $=-\left[\dfrac{X(l_1l_2)_XF_X^2}{2}+\dfrac{Y(l_1l_2)_XF_XF_Y}{2}\right]-\dfrac{X(l_1l_2)_XF_X^2}{2}$

\item $-\left[Xl_1l_2F_{XX}F_X+Yl_1l_2F_{XY}F_X\right]-\dfrac{Y(l_1l_2)_YF_X^2}{2}$

\item $+(X+Y)FF_X +\dfrac{3(l_1l_2)_XF_XF}{2}$

\item $=-2(l_1l_2)_XFF_X-3l_1l_2F_X^2-\dfrac{\left[X(l_1l_2)_XF_X^2+Y(l_1l_2)_YF_X^2\right]}{2}$

\item $+(X+Y)FF_X +\dfrac{3(l_1l_2)_XF_XF}{2}$

\item $-2(l_1l_2)_XFF_X  -3l_1l_2F_X^2- -l_1l_2F_X^2+(X+Y)FF_X +\dfrac{3(l_1l_2)_XFF_X}{2}$

\item $=-l_1l_2F_X^2-3l_1l_2F_X^2-2(l_1l_2)_XFF_X+(X+Y)FF_X+\dfrac{3(l_1l_2)_XFF_X}{2}$

\item $=-4l_1l_2F_X^2$.
\end {description}

It follows finally that

\begin{description}
\item $XGF_X+YGF_Y+F[X(h^{\prime}F_X+k^{\prime}F)-Y(hF_X+kF)]$

\item $=4GF-4l_1l_2F_X^2F=0$ \qed
\end{description}

\vspace{0.2in}
\noindent In particular, the line at infinite is not invariant for $\cal F$. Let us observe that the 16 points of $\{G=0\}\cap S$ are  singularities of $\cal F$ having tangent multiplicity equal to 2. These are {\it all} the singularities of $\cal F$ along $S$. In fact, from \cite{B} we have  that  $(d({\cal F})+2).d(S)= tang(\cal F,S)+S.S$\,\,({\it d} stands for {\it degree}), so that $tang(\cal F,S)=$ 32. But this is already achieved at the singularities of $\{G=0\}\cap S$.

In order to prove that $\widehat S$ is a leaf of $\widehat {\cal F}$ we have to check at each of the points $P_1\dots,P_{16}$ the following conditions (see \cite {Sa}):

\vspace {0.1in}
\begin {description}
\item 1. $A=\dfrac{\partial G}{\partial Y}$
\item 2. $AF_X+BF_Y=0 $
\item 3. $A_Y=\dfrac{G_{YY}}{2}+\dfrac{F_{YY}G_Y}{2F_Y} $
\item 4. $B_Y-A_X=-\dfrac{F_{XY}G_Y}{F_Y}-G_{XY} $
\item 5. $B_X=-\dfrac{G_{XX}}{2}-\dfrac{F_{XX}G_Y}{2F_Y} $
\end{description}

Let us notice that at the points $P_1,\dots,P_4$ we have $(l_1l_2)_X=(l_1l_2)_Y$ and $F_X=F_Y$.

\vspace{0.1 in}
\noindent {\bf Condition 1.}  We have $A=0$ and $\dfrac{\partial G}{\partial Y}=(l_1l_2)_YF_X^2+2l_1l_2F_XF_{XX}=0$ at $P_5,\dots,P_{16}$.
\noindent At $P_1,\dots,P_4$, we see that $A=hF_X=(l_1l_2)_YF_X^2=\dfrac{\partial G}{\partial Y}$

\vspace {0.1in}
\noindent {\bf Condition 2.}  $A=B=0$ at $P_5,\dots,P_{16}$ implies $AF_X+BF_Y=0$ at these points.  As for $P_1,\dots,P_4$:

\vspace {0.1in}
\noindent  $AF_X+BF_Y= hF_X^2-(l_1l_2F_X)_XF_XF_Y=(l_1l_2)_XF_YF_X^2-(l_1l_2)_XF_X^2F_Y=0$

\vspace{0.1in}
\noindent {\bf Condition 3.} We have at $P_5,\dots,P_{16}$: $\dfrac {G_{YY}}{2}=l_1l_2F_{XY}^2$, $\dfrac {F_{YY}G_Y}{2F_Y}=0$ and $A_Y=l_1l_2F_{XY}^2$. At $P_1,\dots,P_4$: 

\vspace {0.1in}
\noindent  $A_Y=\dfrac{(l_1l_2)_{XY}F_YF_X}{2}
+\dfrac{(l_1l_2)_XF_{YY}F_X}{2}
+\dfrac{(l_1l_2)_{YY}F_X^2}{2}+\dfrac{(l_1l_2)_XF_{XY}F_X}{2}+$

\noindent $(l_1l_2)_YF_{XY}F_X+\dfrac{(l_1l_2)_YF_YF_{XY}}{2}+\dfrac{(l_1l_2)_YF_XF_{XY}}{2}-F_XF_Y-\dfrac{(l_1l_2)_XF_{XY}F_Y}{2}=$

\noindent $ F_X^2+\dfrac{(l_1l_2)_XF_{YY}F_X}{2}+2(l_1l_2)_YF_XF_{XY}$ and

\noindent $\dfrac {G_{YY}}{2}+\dfrac{F_{YY}G_Y}{2F_Y}= \dfrac{(l_1l_2)_{YY}F_X^2}{2}+2(l_1l_2)_YF_XF_{XY}+\dfrac{(l_1l_2)_YF_{YY}F_X^2}{2F_Y}=$

\noindent $F_X^2+2(l_1l_2)_YF_XF_{XY}+\dfrac{(l_1l_2)_XF_XF_{YY}}{2}$

\vspace{0.1in}
\noindent {\bf Condition 4.}  At $P_5,\dots,P_{16}$ we have $B_Y-A_X=-2l_1l_2F_{XY}F_{XX}$ and $-\dfrac{F_{XY}G_Y}{F_Y}-G_{XY}=-G_{XY}=-2l_1l_2F_{XY}F_{XX}$. As for $P_1,\dots,P_4:$

\noindent $B_Y-A_X=-(l_1l_2)_{XY}F_X^2-(l_1l_2)_XF_{XY}F_X-(l_1l_2)_YF_XF_{XX}-(l_1l_2)_XF_{XY}F_X+$

\vspace{0.1in} 
\noindent $\dfrac{(l_1l_2)_{XX}F_YF_X}{2}+\dfrac{(l_1l_2)_XF_YF_{XX}}{2}-\dfrac{(l_1l_2)_{XX}F_YF_X}{2}-\dfrac{(l_1l_2)_XF_{XY}F_X}{2}-\dfrac{(l_1l_2)_{XY}F_X^2}{2}-\dfrac{(l_1l_2)_YF_{XX}F_X}{2}-(l_1l_2)_XF_XF_{XY}-\dfrac{(l_1l_2)_XF_YF_{XX}}{2}-\dfrac{(l_1l_2)_YF_XF_{XX}}{2}+F_X^2+\dfrac{(l_1l_2)_XF_{XY}F_X}{2}=-2F_X^2-3(l_1l_2)_XF_XF_{XY}-2(l_1l_2)_YF_XF_{XX}$; this is exactly $\dfrac{F_{XY}G_Y}{F_Y}-G_{XY}$

\vspace{0.1 in}
\noindent {\bf Condition 5.}  At  $P_5,\dots,P_{16}$  we have $B_X=h^{\prime}F_{XX}=-l_1l_2F_{XX}^2  =- \dfrac{G_{XX}}{2}= -\dfrac{G_{XX}}{2}-\dfrac{F_{XX}G_Y}{2F_Y}$

\vspace{0.1 in}
\noindent  As for $P_1,\dots,P_4$:$B_X=-2F_X^2-2(l_1l_2)_XF_XF_{XX}-(l_1l_2)_XF_XF_{XX}+F_X^2+\dfrac{(l_1l_2)_XF_XF_{XX}}{2}=-F_X^2-\dfrac{5(l_1l_2)_XF_XF_{XX}}{2}$.

\vspace{0.1 in}
\noindent On the other hand, $-\dfrac{G_{XX}}{2}-\dfrac{F_{XX}G_Y}{2F_Y}=-F_X^2-2(l_1l_2)_XF_XF_{XX}-\dfrac{(l_1l_2)_YF_XF_{XX}}{2}=-F_X^2-\dfrac{5(l_1l_2)_XF_XF_{XX}}{2}$.

\vspace {0.1in}
\noindent {\bf Other examples}

\vspace{0.1in}
The example above fits into to a family where the points $P_5,\dots,P_{16}$ are the zeros along $\{F=0\}$ of a cubic $C$ which satisfies $C(P_j)= F_X(P_j))$ for $j=1,2,3,4$. We put

\vspace{0.1in}
\begin{description}
\item $G=l_1l_2C^2$
\item $h=\dfrac{(l_1l_2)_YC}{2}+\dfrac{(l_1l_2)_XF_Y}{2}+l_1l_2C_Y=\dfrac{(l_1l_2)_X(C+F_Y)}{2}+l_1l_2C_Y$
\item $k=-C-\dfrac{(l_1l_2)_XC_Y}{2}$
\item $h^{\prime}=-\dfrac{(l_1l_2)_YC}{2}-\dfrac{(l_1l_2)_XF_X}{2}-l_1l_2C_X=-\dfrac{(l_1l_2)_X(C+F_X)}{2}-l_1l_2C_X$
\item $k^{\prime}=C+\dfrac{(l_1l_2)_XC_X}{2}$
\item $A=hF_X+kF$
\item $B=h^{\prime}F_X+k^{\prime}F$
\end{description}

\vspace{0.1in}
We have the same Theorem as before; the computations for the proof are the same. It is clear that if $C(P_j)=\lambda F(P_j)$ for $j=1,2,3,4$ and some $\lambda\ne0$ we may replace $C$ by $\dfrac{C}{\lambda }$.

We do not know if there are cubics (whose zeroes along $\{F=0\}$ are different from $P_1\dots P_4$) for which the Theorem is untrue.

\end {document}